\documentclass[a4paper]{article}

\usepackage[utf8]{inputenc}
\usepackage[english]{babel}
\usepackage[square,numbers]{natbib}

\usepackage{graphicx}
\usepackage{amsthm}

\newtheorem{thm}{Theorem}{}
\newtheorem{lem}{Lemma}{}
\newtheorem{prop}{Proposition}{}
\newtheorem{cor}{Corollary}{}
\newtheorem{rem}{Remark}{}
\newtheorem{defn}{Definition}{}
\newtheorem{alg}{Algorithm}{}
\newtheorem{exmp}{Example}{}

\newenvironment{pf}{\begin{proof}}{\end{proof}}

\usepackage{amsmath, amssymb}

\def\N{\mathbb{N}}
\def\Z{\mathbb{Z}}

\def\R{\mathbb{R}}
\def\C{\mathbb{C}}

\def\rank{\operatorname{rank}}

\def\hmax{h_{\operatorname{max}}}
\def\lmax{\lambda_{\operatorname{max}}}

\newcommand{\norm}[1]{\left\lVert#1\right\rVert}
\newcommand{\abs}[1]{\left\lvert#1\right\rvert}
\newcommand{\mat}[1]{\begin{bmatrix}#1\end{bmatrix}}

\newcommand\Item[1][]{%
  \ifx\relax#1\relax  \item \else \item[#1] \fi
  \abovedisplayskip=0pt\abovedisplayshortskip=0pt~\vspace*{-\baselineskip}}

\title{Singular Value Assignment for Nonuniformly Sampled Systems: Stabilization and Control}
\date{June 2020}
\author{Ufuk Sevim \thanks{Control and Automation Engineering Department, Istanbul Technical University, Istanbul, Turkey } \\ ufuk.sevim@itu.edu.tr \and Leyla Goren-Sumer \footnotemark[1] \\ leyla.goren@itu.edu.tr }

\begin{document}

	
	



	
	\maketitle
	\begin{abstract}
		A state feedback controller design is proposed to guarantee stability of a nonuniformly sampled system for arbitrary selections of sampling periods within an interval, assuming the controller can select the sampling period. It is also shown that such an interval always exists when the continuous time system - that is to be sampled nonuniformly - is stabilizable. The proposed method uses singular value assignment to guarantee the stability. Also, special selections for desired singular values are provided to demonstrate that this method can also be used to improve the performance at the expense of the maximum sampling period.
	\end{abstract}
	
	

\section{Introduction}
The study of nonuniformly sampled systems has been a very active research area in the last decade (\cite{Hetel2017}), mainly due to the practical applications for recently arising topics such as networked and embedded control systems (\cite{Suh2008}, \cite{Baillieul2007}), event based and self-triggered control (\cite{Lunze2010}). These topics find applications in a broad range of areas due to their advantages like flexible architecture, reduced installation and maintenance costs (\cite{Hespanha2007}). Nonuniformly sampled systems can provide a modeling abstraction and framework for analysis and controller design for such systems.

The stability results and control techniques for sampled-data systems are well-established when the sampling interval is uniform, i.e. constant for all times (\cite{Kuo1995}, \cite{Chen1995}). However, these stability conditions are not applicable in the case of nonuniform sampling. It is well-known that arbitrarily varying sampling periods may lead to instability even if the system is stable for each admissible sampling period (\cite{Hetel2017}).

There are many stability conditions and controller design techniques for nonuniformly sampled systems exist in the literature (\cite{Lin2009}, \cite{Fujioka2010}, \cite{Seuret2012}, \cite{Zhang2001}). Many of the works in the literature features some derivation of common quadratic Lyapunov function (CQLF) and solving LMIs, which is very effective but lacks intuition (\cite{Haimovich2013}).

In our work we considered the problem of stabilizing a continuous time system under arbitrarily-sampled state feedback. We provided a method to find a coordinate transformation which allows us to assign singular values of the discretized closed loop system less than 1. This guarantees stability since the length of the state vector always decreases in this particular coordinate system for arbitrary sampling period. We showed that such a coordinate transformation always exists if the continuous system is stabilizable. We also showed that this is equivalent to finding a CQLF for the discretized system matrices. 




The proposed method is based on the main theorem given in \cite{Martin2009}, which gives the necessary and sufficient conditions for the existence of a state feedback gain to assign the desired singular values of the closed-loop system. The authors also provide an algorithm to calculate such a state feedback gain.

We also introduced particular singular value selections to improve the performance of the closed loop system at the expense of the maximum sampling time. 

To the best of our knowledge, this is the first work that uses singular value assignment to stabilize and control non-uniformly sampled systems.

\subsection{Notation and Terminology}
For $x \in \C^n$, $\abs{x}$ denotes any vector norm and for $A \in \C^{n \times n}$, $\norm{A}$ denotes a norm induced by some vector norm $\abs{\cdot}$, i.e. $\norm{A}:=\max_{x \in \C^n} \abs{Ax}/\abs{x}$. $\bar{\sigma}(A)$ denotes the maximum singular value of $A$.

\section{Problem Formulation}
Consider the continuous time system
\begin{equation}
	\label{eq:cont_ol}
	\dot{x}(t) = Ax(t) + Bu(t)
\end{equation}
where $x(t) \in \R^n$, $u(t) \in \R^m$, $A \in \R^{n \times n}$, $B \in \R^{n \times m}$, $\rank B = m \leq n$, and $(A, B)$ is stabilizable.

Define the sequence of time instances $\{t_k\}_{k \in \N}$ where
\[ 0 = t_0 < t_1 < \dots < t_k < \dots \]
with $\lim_{k \to \infty} t_k = \infty$. We assume that $\{t_k\}_{k \in \N}$ are known a priori and sampling intervals $h_k := t_{k+1} - t_k$ are bounded, i.e. $h_k \in (0, \hmax), \forall k \in \N$ for some $\hmax > 0$.

The objective is to design a sampled state feedback controller $K(t_k) \in \R^{m \times n}$ such that the system
\begin{equation}
	\label{eq:sd_cl}
	\dot{x}(t) = A x(t) + B K(t_k) x(t_k), ~~ \forall t \in [t_k, t_{k+1})
\end{equation}
is stable for any selection of $\{t_k\}_{k \in \N}$ as long as $h_k \in (0, \hmax)$.

Consider the discretized model of \eqref{eq:cont_ol}
\begin{equation}
	\label{eq:discrete_ol}
	x_{k+1} = F_k x_k + G_k u_k
\end{equation}
where $u_k := u(t_k)$, 
\[ F_k := e^{A h_k} ~~ \text{and} ~~ G_k := \left( \int_0^{h_k} e^{A \tau} d \tau \right) B. \]
The following result is standard:
\begin{lem}
	\label{lem:discretization}
	Define $K_k := K(t_k)$ where $K(t_k)$ is given in \eqref{eq:sd_cl}. Then the solutions of \eqref{eq:sd_cl} and
	\begin{equation}
		\label{eq:discrete_cl}
		x_{k+1} = (F_k + G_k K_k) x_k
	\end{equation}
	are same at the sampling instances, i.e. $x_k = x(t_k), \forall k \in \N$, assuming $x(t_0) = x_0$.
\end{lem}

\begin{pf}
	We prove with induction. $x(t_0) = x_0$ by assumption. Assume that $x(t_k) = x_k$, so
	\begin{align*}
		x(t_{k+1}) &= e^{A (t_{k+1} - t_k)} x(t_k) + \int_{t_k}^{t_{k+1}} e^{A(t_{k+1} - \eta)} B u(\eta) d \eta \\
		&= e^{A h_k} x(t_k) + \left( \int_0^{h_k} e^{A \tau} d \tau \right) B u(t_k) \\
		&= F_k x_k + G_k K_k x_k \\
		&= x_{k+1}
	\end{align*}
\end{pf}

So the problem becomes finding $K_k$ and $\hmax$ such that system \eqref{eq:discrete_cl} is guaranteed to be stable for $h_k \in (0, \hmax)$. We are going to use singular value assignment to guarantee the stability of \eqref{eq:discrete_cl}.

\section{Background}
It is well-known that the stability conditions for nonuniformly sampled systems are not trivial (\cite{Liberzon1999}). There are many sufficient stability results exist, such as \cite{Lin2009}, \cite{Kao2013}, \cite{Seuret2012}, \cite{Liberzon1999a}.

We use the following definition and sufficient condition for stability.
\begin{defn}
	The system
	\begin{equation}
		\label{eq:sd_ol}
		x_{k+1} = F_k x_k
	\end{equation}
	is globally asymptotically stable if, for any $x_0 \in \R^n$
	\begin{equation*}
		\lim_{k \to \infty} \abs{x_k} = 0
	\end{equation*}
	for some vector norm $\abs{\cdot}$.
\end{defn}

\begin{thm}
	\label{thm:stability_norm}
	If there exists a matrix norm $\norm{\cdot}$ induced over some vector norm $\abs{\cdot}$ such that $\norm{F_k} < 1, \forall k \in \N$, then the system \eqref{eq:sd_ol} is globally asymptotically stable.
\end{thm}

\begin{pf}
	Let $V(x_k) := \abs{x_k}$. Then,
	\[ V(x_{k+1}) = \abs{x_{k+1}} \leq \norm{F_k} \abs{x_k} < \abs{x_k} = V(x_k). \]
	Hence $V(\cdot)$ is a Lyapunov function for the system \eqref{eq:sd_ol}.
\end{pf}

\begin{lem}
	\label{lem:norm_similarity}
	Let $\norm{\cdot}$ be a matrix norm induced over some vector norm $\abs{\cdot}$. Then for any $A \in \R^{n \times n}$
	\begin{equation}
		\norm{A}_T := \norm{T^{-1} A T}
	\end{equation}
	is an induced norm over the vector norm
	\begin{equation}
		\abs{x}_T := \abs{T^{-1} x}
	\end{equation}
	where $x \in \R^n$ and $T \in \R^{n \times n}$ is an invertible matrix. 
\end{lem}

\begin{pf}
	It is easy to see that $\abs{\cdot}_T$ is a vector norm. Then,
	\[
		\sup_{x \in \R^n} \frac{\abs{A x}_T}{\abs{x}_T} = \sup_{y \in \R^n} \frac{\abs{A T y}_T}{\abs{T y}_T} = \sup_{y \in \R^n} \frac{\abs{T^{-1} A T y}}{\abs{T^{-1} T y}} = \norm{A}_T
	\]
\end{pf}

\begin{cor}
	\label{cor:stability_sv}
	If there exists an invertible matrix $T \in \R^{n \times n}$ such that
	\begin{equation}
		\label{ineq:singular_value}
		\bar{\sigma}(T^{-1} F_k T) < 1, ~~ \forall k \in \N,
	\end{equation}
	then the system \eqref{eq:sd_ol} is globally asymptotically stable, where $\bar{\sigma}(\cdot)$ is the maximum singular value.
\end{cor}

\begin{pf}
	It is well-known that the maximum singular value of a matrix is the induced matrix norm over Euclidean vector norm. Then the result follows from Lemma \ref{lem:norm_similarity}.
\end{pf}

The inequality \eqref{ineq:singular_value} is actually equivalent to the common quadratic Lyapunov function condition, as shown in the following proposition.
\begin{prop}
	There exists a $P > 0$ such that
	\begin{equation}
		F_k^T P F_k - P < 0, ~~ \forall k \in \N
	\end{equation}
	if and only if there exists an invertible $T \in \R^{n \times n}$ such that
	\[ \bar{\sigma}(T^{-1} F_k T) < 1, ~~ \forall k \in \N. \]
\end{prop}
\begin{pf}
	($\Rightarrow$) Since $P > 0$, there exists invertible $U \in \R^{n \times n}$ such that $P = U^T U$. So,
	\begin{align*}
		F_k^T U^T U F_k - U^T U &< 0 \\
		U^{-T} F_k^T U^T U F_k U^{-1} - I &< 0 \\
		\lmax(U^{-T} F_k^T U^T U F_k U^{-1}) &< 1 \\
		\bar{\sigma} (U F_k U^{-1}) &< 1
	\end{align*}
	Now, take $T = U^{-1}$. \\
	($\Leftarrow$) Similarly, go backwards and take $P = T^{-T} T^{-1}$.
\end{pf}	

\cite{Martin2009} gives the necessary and sufficient conditions for assigning singular values using state feedback gain, as well as an algorithm for finding such a gain. The following theorem is given without proof for completeness.

\begin{thm}
	\label{thm:singular_value_assignment}
	For any given $n \times n$ matrix $A$ and $n \times m$ full rank matrix $B$ [over reals] with $1 \leq m \leq n$, let the singular values of
	\begin{equation}
		(I - B(B^T B)^{-1} B^T)A
	\end{equation}
	be
	\[ 0 = a_1 = \dots = a_m \leq a_{m+1} \leq \dots \leq a_n. \]
	For any given set of values
	\[ 0 \leq s_1 \leq s_2 \leq \dots \leq s_n \]
	there exists a real $K$ such that the singular values of $A+BK$ are $\{s_1, \dots, s_n\}$ if and only if
	\begin{equation}
		a_j \leq s_j \leq a_{j+m}
	\end{equation}
	for all $j=1,\dots,n$ (with the convention that $a_j = \infty$ if $j > n$).
\end{thm}

It is important to note that if such a $K$ exists it is not unique. The algorithm given in \cite{Martin2009} calculates only one such $K$, except for the special case of $\rank B = 1$ and $a_1, \dots, a_n$ are distinct, in which case there are finitely many $K$.

The following well-known results are given for completeness.
\begin{lem}
	\label{lem:symmetric_sum}
	For symmetric real matrices $P$ and $Q$, the following inequality holds:
	\begin{equation}
		\lmax(P+Q) \leq \lmax(P) + \lmax(Q)
	\end{equation}
\end{lem}
\begin{pf}
	For symmetric real matrices we can write:
	\begin{align*}
		\lmax(P+Q) &= \sup_{\abs{v}=1} v^T (P + Q) v \\
				   &\leq \sup_{\abs{v}=1} v^T P v + \sup_{\abs{v}=1} v^T Q v \\
				   &= \lmax(P) + \lmax(Q)
	\end{align*}
\end{pf}


\begin{lem}
	\label{lem:tridiag}
	The tridiagonal symmetric $n \times n$ matrix
	\[ M = \mat{0 & 1 & & & \\ 1 & 0 & 1 & & \\ & \ddots & \ddots & \ddots & \\ & & 1 & 0 & 1 \\ & & & 1 & 0} \]
	has eigenvalues $\lambda_i = 2 \cos(i \pi / (n+1)), ~ i = 1,2,\dots,n$.
\end{lem}

\begin{pf}
	Let $M x = \lambda x$ where $x^T := \mat{x_1 & \dots & x_n}$. Then,
	\begin{equation}
		\label{eq:tridiag}
		x_{k+2} = \lambda x_{k+1} - x_k, ~~~~ k=0,1,\dots,n-1
	\end{equation}
	with the convention $x_0 = x_{n+1} = 0$. It is easy to show that the trigonometric identity
	\[ \sin[(k+2) u] = 2 \cos (u) \sin [(k+1) u] - \sin (ku) \]
	holds for any $u \in \R$. Thus, \eqref{eq:tridiag} is satisfied by selecting
	\[ x_k := \sin (ku) ~~~~ \text{and} ~~~~ \lambda := 2 \cos (u). \]
	The condition $x_{n+1} = \sin [(n+1) u] = 0$ can be satisfied by selecting
	\[ u := i \pi / (n+1), ~~~~ i = 1,2,\dots,n \]
	which proves the claim.
\end{pf}

\section{Stabilization}
The following result is immediate from Theorem \ref{thm:singular_value_assignment} and Corollary \ref{cor:stability_sv}.

\begin{cor}
	\label{cor:main}
	If there exists an invertible matrix $T \in \R^{n \times n}$ such that
	\begin{equation}
		\label{eq:sva_main}
		\bar{\sigma}\left( \left( I - \hat{G}(h) \left[ \hat{G}^T(h) \hat{G}(h) \right]^{-1} \hat{G}^T(h) \right) \hat{F}(h) \right) < 1
	\end{equation}
	for all $h \in (0, \hmax)$ for some $\hmax > 0$, where
	\[ \hat{F}(h) := T^{-1} e^{Ah} T ~~ \text{and} ~~ \hat{G}(h) := T^{-1} \left( \int_{0}^{h} e^{A \tau} d \tau \right) B, \]
	then there exists $\hat{K}(h) \in \R^{m \times n}$ such that
	\[ \bar{\sigma} \left( \hat{F}(h) + \hat{G}(h) \hat{K}(h) \right) < 1 \]
	for all $h \in (0, \hmax)$. Also the state feedback gain
	\begin{equation}
		K_k := \hat{K}(h_k) T^{-1}
	\end{equation}
	guarantees the stability of \eqref{eq:discrete_cl} for any given $\{t_k\}_{k \in \N}$ such that $h_k \in (0, \hmax)$.
\end{cor}

Corollary \ref{cor:main} assumes that $\hat{G}(h)$ is full-rank for all $h \in (0, \hmax)$. The following lemma gives necessary and sufficient conditions for $\hat{G}(h)$ to be a full-rank matrix.
\begin{lem}
	\label{lem:full_rank}
	$\hat{G}(h)$ is full-rank for all $h \in (0, \hmax)$ if and only if $B$ is full-rank and
	\[ \lambda h \neq 2k\pi j, ~~ \forall k \in \Z \backslash \{0\}, ~~ \forall h \in (0, \hmax) \]
	where $\lambda$ is any eigenvalue of $A$ and $j = \sqrt{-1}$.
\end{lem}

\begin{pf}
	We need to show that $\int_{0}^{h} e^{A \tau} d \tau$ is singular if and only if $\lambda h = 2k\pi j$. Since
	\[ \int_{0}^{h} e^{A \tau} d \tau = \sum_{n=1}^\infty \frac{A^{n-1} h^n}{n!} \]
	its eigenvalues are in the form
	\[ \phi(h) := \begin{cases}
			h & \lambda = 0 \\
			\frac{e^{\lambda h} - 1}{\lambda} & \lambda \neq 0
		\end{cases}
	\]
	as a consequence of Spectral Mapping Theorem. So, $\phi(h)=0$ if and only if $e^{\lambda h} = 1$, which concludes the result.
\end{pf}

Lemma \ref{lem:full_rank} is a special case of the well-known pathological sampling phenomenon. A sampling period $h$ is said to be pathological if there exists $k \in \Z \backslash \{0\}$ such that
\[ h(\lambda - \mu) = 2k\pi j \]
where $\lambda$ and $\mu$ are any pair of eigenvalues of $A$.

\begin{lem}
	Let $(A,B)$ be controllable. If $h$ is not pathological, then $\left( F(h), G(h) \right)$ is controllable. Furthermore, if the system has single input, i.e. $\rank B = 1$, and $\left( F(h), G(h) \right)$ is controllable, then $h$ is not pathological.
\end{lem}
\begin{pf}
	See \cite{Chen1984} Appendix D.
\end{pf}

When a pathological sampling occurs, there is a possibility for the discrete time system to lose its controllability, i.e. some eigenvalues of $F(h)$ cannot be moved by the state feedback. It is important to note that losing controllability does not induce any other conditions for singular value assignment other than Lemma \ref{lem:full_rank} and Theorem \ref{thm:singular_value_assignment}.

However, it is well-known that the maximum singular value cannot be less than the spectral radius for any matrix. Therefore, if an uncontrollable eigenvalue of $F(h)$ has a magnitude greater than 1, then the condition \eqref{eq:sva_main} is not satisfied for any similarity transformation matrix, which means $\hmax$ must be less than the minimum pathological sampling period due to the continuity of singular values. On the other hand, if all uncontrollable eigenvalues of $F(h)$ lies in the unit disk, it can still be possible to assign singular values to less than 1.

From these discussions we can see that no special treatment is required for the case of pathological sampling as long as $\hat{G}(h)$ is full-rank. Now, we present the main lemma of this paper.

\begin{lem}
	\label{lem:main}
	Let $(A, B)$ be given in \eqref{eq:cont_ol} and there exists $K \in \R^{m \times n}$ such that $A+BK$ is stable. Let
	\[ S := T^{-1} (A + BK) T \]
	and $S+S^T<0$. Then, there exists an interval $(0,\hmax)$ such that \eqref{eq:sva_main} is satisfied with this $T$.
\end{lem}

\begin{pf}
	Using the property
	\[ e^{Ah} = I + \left( \int_{0}^{h} e^{A \tau} d \tau \right) A \]
	and by substituting $T^{-1} A T = S - T^{-1} B K T$ we obtain
	\begin{equation}
		\hat{F}(h) = I + T^{-1} \left( \int_{0}^{h} e^{A \tau} d \tau \right) T S - \hat{G}(h) K T
	\end{equation}
	Then,
	\begin{align*}
		[I - \hat{G} (\hat{G}^T \hat{G})^{-1} \hat{G}^T] \hat{F} = [I - \hat{G} (\hat{G}^T \hat{G})^{-1} \hat{G}^T] \left[ I + T^{-1} \left( \int_{0}^{h} e^{A \tau} d \tau \right) T S \right]
	\end{align*}
	The left matrix is an orthogonal projection matrix, so it always has a maximum singular value of 1. Therefore it is sufficient to show that
	\begin{equation}
	\label{eq:msv_proj}
	\bar{\sigma}\left[ I + T^{-1} \left( \int_{0}^{h} e^{A \tau} d \tau \right) T S \right] < 1
	\end{equation}
	holds for a small enough $h$. Using the series expansion of the above matrix and the definition of maximum singular value, the condition above is equivalent to
	\begin{equation}
	\lmax \left[ I + h (S^T + S) + h^2 Q(h) \right] < 1
	\end{equation}
	where $Q(h)$ is the remaining part, which is symmetric and entire, i.e. it is analytic everywhere. Using Lemma \ref{lem:symmetric_sum} and Spectral Mapping Theorem, we can see that
	\begin{align*}
		\lmax \left[ I + h (S^T + S) + h^2 Q(h) \right] &\leq \lmax \left[ I + h (S^T + S) \right] + h^2 \lmax \left[ Q(h) \right] \\
		&= 1 + h \lmax (S^T + S) + h^2 \lmax \left[ Q(h) \right] \\
		&< 1
	\end{align*}
	is sufficient for the proof, which is equivalent to
	\[ \lmax (S^T + S) + h \lmax \left[ Q(h) \right] < 0 \]
	Select any $\bar{h} > 0$. Since $Q(h)$ is entire, there exists $M > 0$ such that $\lmax \left[ Q(h) \right] < M$ for all $h \in (0, \bar{h})$. Finally, select
	\[ \hmax = \min \left( \bar{h}, \frac{-\lmax (S^T + S)}{M} \right) \]
	Since $S+S^T$ has negative eigenvalues, this selection guarantees \eqref{eq:sva_main} holds for all $h \in (0, \hmax)$, which proves the claim.
\end{pf}

As a consequence of the above Lemma, the main theorem of the paper can be given.

\begin{thm}
	\label{thm:main}
	Let $(A, B)$ be stabilizable. Then there exists $\hmax > 0$ and $K(t_k) \in \R^{m \times n}$ such that the system
	\[ \dot{x}(t) = A x(t) + B K(t_k) x(t_k), ~~ \forall t \in [t_k, t_{k+1}) \]
	is stable for any selection of $\{t_k\}_{k \in \N}$ as long as $h_k \in (0, \hmax)$.
\end{thm}
\begin{pf}
	Let $K$ be a stabilizing controller for the continuous system, i.e. $A+BK$ is Hurwitz. Then, there exists $P > 0$ such that
	\[ (A+BK)^T P + P (A+BK) < 0 \]
	Since $P>0$, there exists an invertible $T \in \R^{n \times n}$ such that $P = T^{-T} T^{-1}$, which yields
	\[ (T^{-1} (A+BK) T)^T + T^{-1} (A+BK) T < 0 \]
	So, using Lemma \ref{lem:main}, Corollary \ref{cor:stability_sv} and Lemma \ref{lem:discretization} we can conclude the result.
\end{pf}

Note that there are infinitely many choices of $T$. We give a few well-known methods for choosing $T$ below.

\begin{rem}
	Both $K$ and $T$ can also be found by solving the Riccati equation. It is well-known that if $(A,B)$ is stabilizable and $(A,Q)$ is observable, then the Riccati equation
	\[ A^T P + P A - P B R^{-1} B^T P + Q = 0 \]
	has a unique stabilizing solution $P > 0$ and the state feedback gain is $K = -R^{-1} B^T P$, where $Q \geq 0$ and $R > 0$ are gains of the quadratic cost function. Since,
	\[ (A+BK)^T P + P (A+BK) = - P B R^{-1} B^T P - Q < 0 \]
	the stability can be guaranteed by selecting $T$ such that $P = T^{-T} T^{-1}$.
\end{rem}

\begin{rem}
	If $K$ can be selected such that $A+BK$ is diagonalizable, e.g. by assigning distinct eigenvalues, $T$ can be selected such that
	$J := T^{-1} (A+BK) T$ has real block diagonal form, where $i$th block $J_i=\lambda_i$ is a scalar if the corresponding eigenvalue is real and
	\[ J_i = \mat{a_i & b_i \\ -b_i & a_i} \]
	if the corresponding eigenvalues are complex conjugate pairs, that is $\lambda_i = a_i \mp j b_i$. So, $J + J^T$ becomes diagonal with negative real eigenvalues, which satisfies Lemma \ref{lem:main}.
\end{rem}

\begin{rem}
	If $A+BK$ is defective, e.g. there are uncontrollable defective eigenvalues, then there exists a real $T$ such that $J$ has real block diagonal form and each block has the Jordan normal form, that is
	\[ J_i = \mat{\lambda_i & 1 & & \\ & \ddots & \ddots & \\ & & \lambda_i & 1 \\ & & & \lambda_i} \]
	if $\lambda_i$ is real and
	\[ J_i = \mat{C_i & I & & \\ & \ddots & \ddots & \\ & & C_i & I \\ & & & C_i} \]
	if $\lambda_i = a_i \mp j b_i$ are complex conjugate pairs where
	\[ C_i = \mat{a_i & b_i \\ -b_i & a_i} \]
	In this case, $J + J^T$ is negative definite if
	\begin{equation}
		\label{cond:defective}
		\operatorname{Re} \lambda_i < -\cos(\pi/(k_i+1))
	\end{equation}
	where $k_i$ is the multiplicity of $\lambda_i$ within the block $J_i$. This fact is immediate since
	\[ J_i + J_i^T = 2 \lambda_i I + M \]
	where $M$ is defined as in Lemma \ref{lem:tridiag} for the real case and it is easy to show for the complex conjugate case.
\end{rem}

\begin{rem}
	Once a suitable $T$ is found, $TV$ is also a suitable similarity transformation matrix where $V$ is an orthogonal matrix, i.e. $V^TV=I$. This can immediately be seen from \eqref{eq:sva_main} as the matrix becomes
	\begin{equation}
	V^T \left( I - \hat{G}(h) \left[ \hat{G}^T(h) \hat{G}(h) \right]^{-1} \hat{G}^T(h) \right) \hat{F}(h) V
	\end{equation}
	which has the same singular values. This means selection of $V$ does not affect stability or assignable singular value regions. 
\end{rem}

An algorithm for finding a stabilizing controller can be given as follows to summarize the process:

\begin{alg}
	The following design algorithm for a stabilizing controller is proposed based on Theorem \ref{thm:main}.
	
	\begin{enumerate}
		\item[1.] Find $K \in \R^{m \times n}$ and $T \in \R^{n \times n}$ with one of the following methods:
		\begin{enumerate}
			\item[a.] Find $K$ such that $A+BK$ is Hurwitz, then solve the Lyapunov equation
			\[ (A+BK)^T P + P (A+BK) + Q = 0 \]
			for some $Q > 0$. Find $T$ such that $P = T^{-T} T^{-1}$.
			\item[b.] Solve the Riccati equation
			\[ A^T P + P A - P B R^{-1} B^T P + Q = 0 \]
			for some $R > 0$ and $Q \geq 0$. Calculate $K = -R^{-1} B^T P$ and $T$ such that $P = T^{-T} T^{-1}$.
			\item[c.] Find $K$ such that $A+BK$ is Hurwitz and calculate $T$ such that $T^{-1} (A+BK) T$ has real Jordan block diagonal form. Check if the condition \eqref{cond:defective} holds on defective case.
		\end{enumerate}
	
		\item[2.] Calculate
		\[ \hat{F}(h) := T^{-1} e^{Ah} T ~~ \text{and} ~~ \hat{G}(h) := T^{-1} \left( \int_{0}^{h} e^{A \tau} d \tau \right) B. \]
		
		\item[3.] Calculate the singular values of
		\[\hat{P}(h) := \left( I - \hat{G}(h) \left[ \hat{G}^T(h) \hat{G}(h) \right]^{-1} \hat{G}^T(h) \right) \hat{F}(h). \]
		Due to Theorem \ref{thm:main}, there exists $\hmax$ such that $\bar{\sigma}(\hat{P}(h)) < 1$ for all $h \in (0, \hmax)$.
		
		\item[5.] Determine the desired singular values according to Theorem \ref{thm:singular_value_assignment}.
		
		\item[6.] Calculate $\hat{K}(h)$ such that $\hat{F}(h) + \hat{G}(h) \hat{K}(h)$ has desired singular values using the algorithm given in \cite{Martin2009}.
		
		\item[7.] The controller $K_k := \hat{K}(h_k) T^{-1}$ is a stabilizing controller for the system \eqref{eq:discrete_cl}.
	\end{enumerate}
\end{alg}

\begin{exmp}
	We consider the example given in \cite{Haimovich2013} with an additional unstable eigenvalue. This example is also given in our previous work \cite{Sevim2016}.
	\[ A = \mat{1 & -2 & 0 \\ 2 & 1 & 0 \\ 0 & 0 & 0.5}, ~~ B = \mat{0.5 \\ 2 \\ 1} \]
	We selected the following state feedback gain to stabilize the system:
	\[ K = \mat{\frac{1128}{289} & \frac{-1064}{289} & \frac{-105}{34}} \]
	
	This controller assigns the closed loop system eigenvalues to $\{-1,-2,-3\}$. We selected the similarity transformation matrix $T$ that diagonalizes $A+BK$.
	
	The nonzero singular values of $\hat{P}(h)$ are calculated as in Figure \ref{fig:ex1_sv_ph}. From the figure, we can conclude that $\hmax \approx 0.62$ and the singular values of $\hat{F}(h) + \hat{G}(h) \hat{K}(h)$ can be arbitrarily assigned anywhere between the solid curves. We selected the dashed lines to assign the singular values, which guarantees the stability of the closed loop system for any selection of sampling time between $(0,\hmax)$.
	
	In Figure \ref{fig:ex1_cl_resp}, the state response of the closed loop system \eqref{eq:sd_cl} under 100 random sampling is given.
	
	\begin{figure}[!ht]
		\centering
		\includegraphics[width=.8\columnwidth]{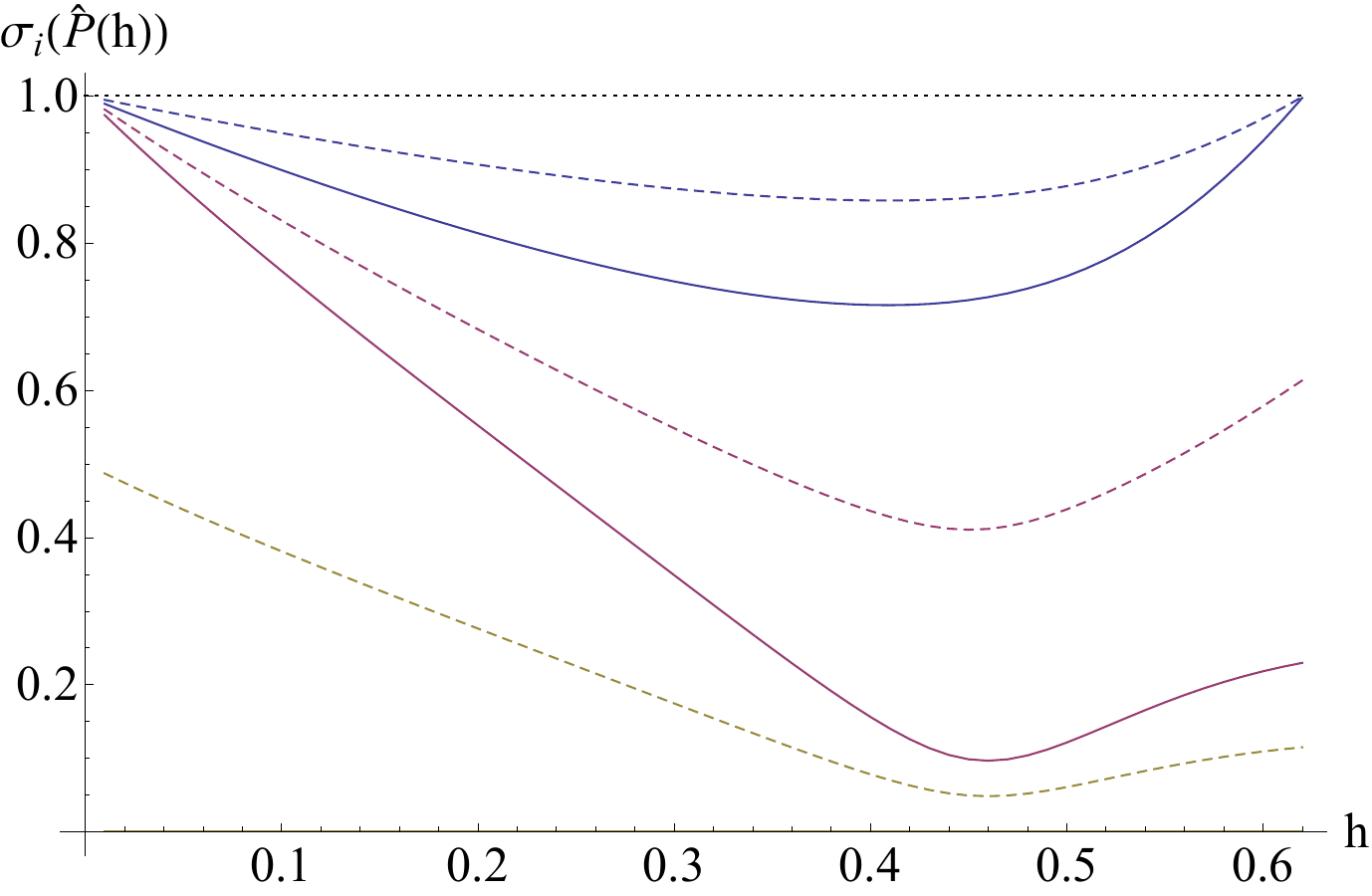}
		\caption{The nonzero singular values of $\hat{P}(h)$ (solid curves) and the desired singular values (dashed curves).}
		\label{fig:ex1_sv_ph}
	\end{figure}
	
	\begin{figure}[!ht]
		\centering
		\includegraphics[width=1\columnwidth]{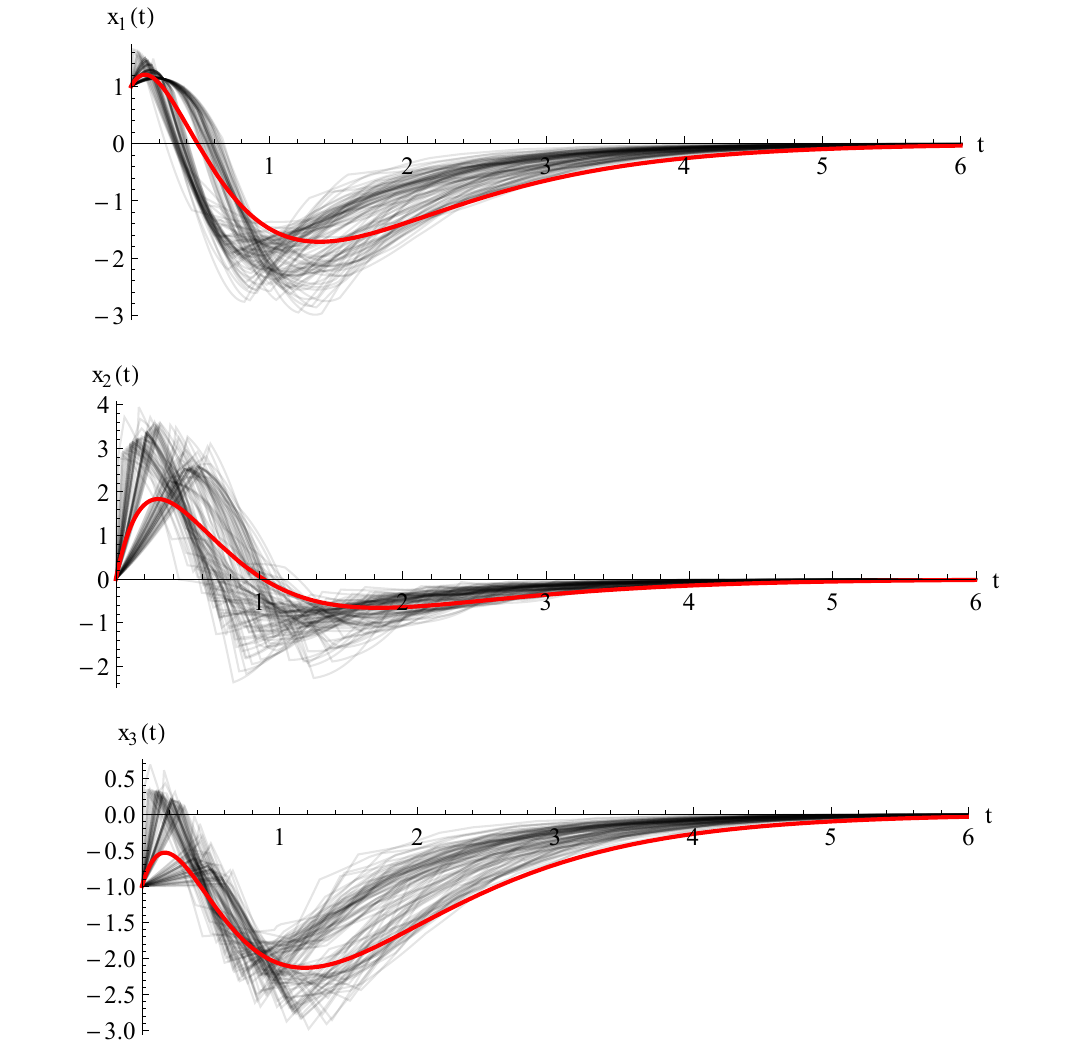}
		\caption{The state response of the closed loop system \eqref{eq:sd_cl} under 100 random sampling where $h_k \in (0, 0.62)$ (grey lines) and the state response of the continuous time system $\dot{x}(t) = (A+BK)x(t)$ (red line).}
		\label{fig:ex1_cl_resp}
	\end{figure}
\end{exmp}

\section{Control}
In this section we proposed a selection of the desired singular values which does not only stabilize the system but also provide a good performance. Here performance is defined by the "average distance" between the continuous and sampled-data system state responses. More precisely, the performance of a controller $K(h)$ can be defined as
\[
	P(K(h)) := E \left[\norm{e(t)}\right] = E \left[ \left( \int_{0}^{\infty} e^T(t) e(t) dt \right) ^ {1/2} \right]
\]
where $E$ denotes the expected value when $\{h_k\}$ are randomly selected, $e(t) := x(t) - x_d(t)$, $x(t)$ is the solution to \eqref{eq:sd_cl}, $x_d(t)$ is the solution to \eqref{eq:cont_ol} with $u(t) = K x_d(t)$ and $x(0) = x_d(0)$.

Assume that $A+BK$ has real distinct eigenvalues. We observed that selecting diagonalizing $T$, i.e. $S$ and $e^{Sh}$ are diagonal, gives good performance if the singular values are assigned to $\sigma_i(e^{Sh})$ provided that the assignment preserves the structure of the closed loop system, namely if $\hat{F}(h) + \hat{G}(h) \hat{K}(h)$ is also diagonally dominant. Since $\hat{K}(h)$ is not unique, the closed loop system may not always become diagonally dominant with the $\hat{K}(h)$ produced by the singular value assignment algorithm. In this case, one has to choose the state feedback which preserves diagonal dominance.

%


\begin{exmp}
	We consider the example given in \cite{Fujioka2010}
	\[ A = \mat{0 & 1 \\ 0 & -0.1}, ~~ B = \mat{0 \\ 0.1} \]
	We selected the following state feedback gain to stabilize the system:
	\[ K = \mat{-20 & -29} \]
	
	This controller assigns the closed loop system eigenvalues to $\{-1,-2\}$. We selected the similarity transformation matrix $T$ that diagonalizes $A+BK$.
	
	The nonzero singular values of $\hat{P}(h)$ are calculated as in Figure \ref{fig:ex2_sv_ph}. From the figure, we can conclude that there is no limit to $\hmax$ for stabilization. However, $\hmax \approx 1.65$ for the ability to assign the singular values of $\hat{F}(h) + \hat{G}(h) \hat{K}(h)$ to $\sigma_i(e^{Sh})$, which is depicted as dashed lines.
	
	In Figure \ref{fig:ex2_cl_resp2}, the state response of the closed loop system \eqref{eq:sd_cl} under 100 random sampling is given. 
%
	
	\begin{figure}[!ht]
		\centering
		\includegraphics[width=.8\columnwidth]{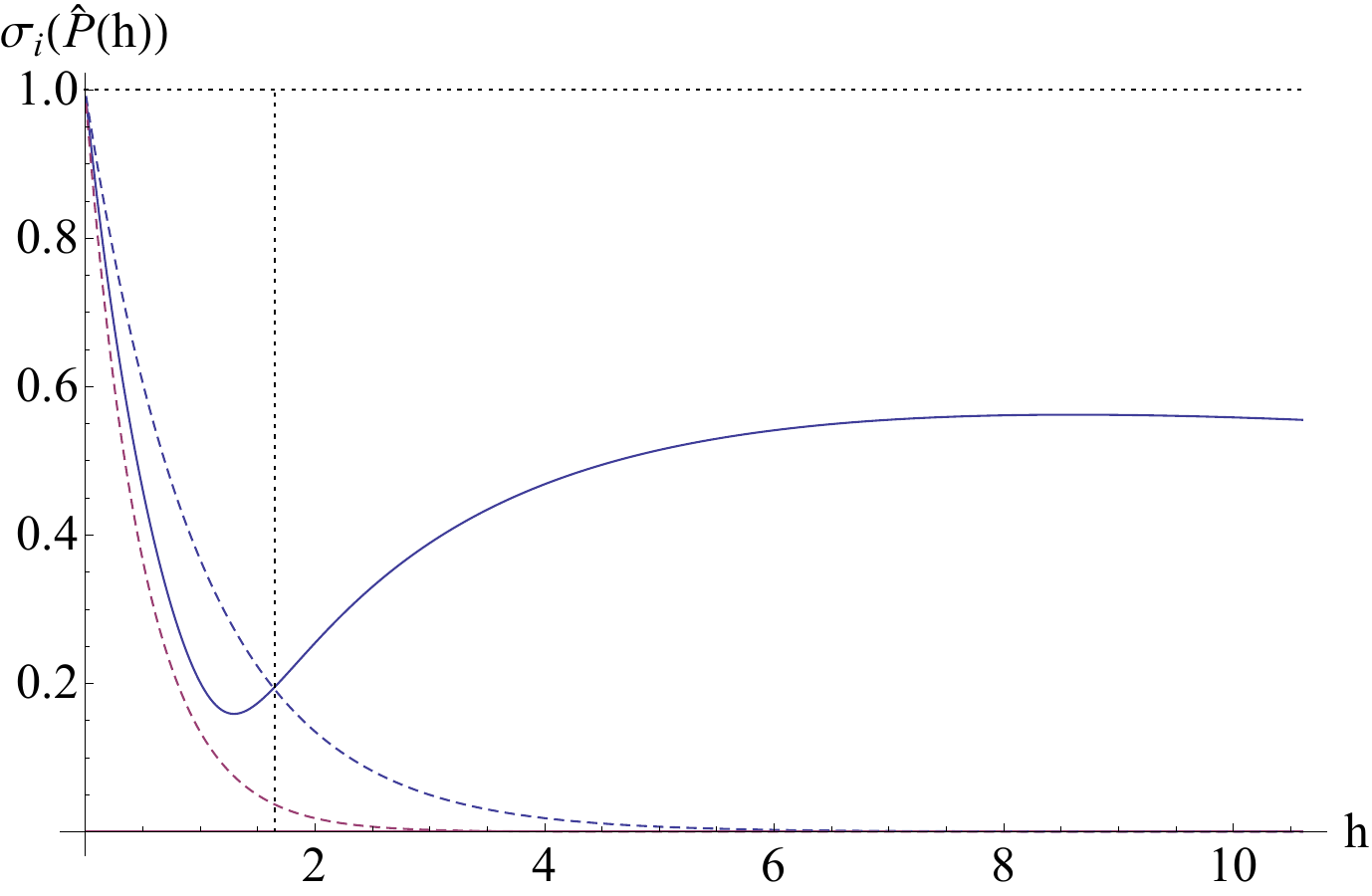}
		\caption{The nonzero singular values of $\hat{P}(h)$ (solid curves) and the desired singular values (dashed curves).}
		\label{fig:ex2_sv_ph}
	\end{figure}
	

	\begin{figure}[!ht]
		\centering
		\includegraphics[width=\columnwidth]{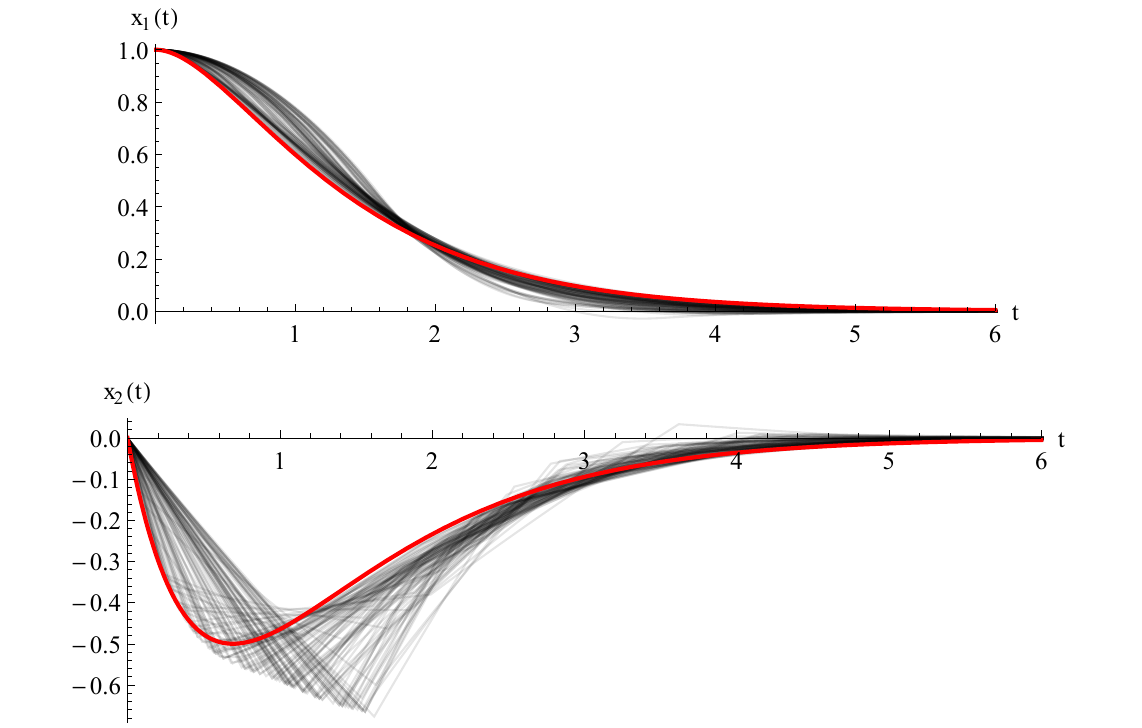}
		\caption{The state response of the closed loop system \eqref{eq:sd_cl} under 100 random sampling where $h_k \in (0, 1.65)$ (grey lines) and the state response of the continuous time system $\dot{x}(t) = (A+BK)x(t)$ (red line).}
		\label{fig:ex2_cl_resp2}
	\end{figure}

	
\end{exmp}

\section{Conclusions}
We considered the stabilization problem of nonuniformly sampled systems with the assumption that sampling periods are known to the controller. We presented a theorem for the existence of a state feedback controller that assigns the singular values of the closed loop system matrix to desired values under a similarity transformation. We provided 3 methods to calculate such a similarity transformation matrix. Also, we provided an algorithm to calculate the sampling-period-varying state feedback controller that guarantees stability for bounded and arbitrarily varying sampling intervals. We illustrated the accuracy of the method in numerical examples. Lastly, we proposed specific singular values to assign to improve the performance of the system.

		
	\bibliographystyle{plain}
	\bibliography{references}
\end{document}